\documentclass[12pt]{article}

\usepackage{amsmath}
\usepackage{amssymb}
\usepackage{amsthm}

\usepackage[dvips]{graphicx}
\usepackage{color}
\usepackage{epsfig}
\usepackage[mathscr]{eucal}

\newtheorem{theorem}{Theorem}

\newtheorem{lemma}[theorem]{Lemma}
\newtheorem{cor}[theorem]{Corollary}

\newcommand{\floor}[1]{\ensuremath{\left \lfloor {#1} \right \rfloor}}

\newcommand{\bv}{{\bf v}}
\newcommand{\bw}{{\bf w}}
\newcommand{\bx}{{\bf x}}
\newcommand{\bP}{{\bf P}}
\newcommand{\bE}{{\bf E}}

\newcommand{\cC}{{\cal C}}
\newcommand{\cE}{{\cal E}}

\newcommand{\cH}{{\cal H}}

\newcommand{\cV}{{\cal V}}
\newcommand{\cX}{{\cal X}}

\newcommand{\ignore}[1]{}

\allowdisplaybreaks

\title{De Bruijn Covering Codes for Rooted Hypergraphs}
\author{Fan Chung \\ \small Department of Mathematics \\ \small University of California, San Diego \\[.2in]
Joshua N. Cooper \\ \small Courant Institute of Mathematics \\ \small New York University}

\begin{document}

\maketitle

\begin{abstract}
What is the length of the shortest sequence $S$ of reals so that the set of 
consecutive $n$-words in $S$ form a covering code for permutations on
$\{1,2, \ldots, n\}$  of radius $R$ ?  
(The distance between two $n$-words is the number 
of transpositions needed to have the same order type.) 
The above problem can be viewed as a special case of
finding a De Bruijn covering
code for a rooted hypergraph. Each edge of a rooted hypergraph contains
a special vertex, called the {\it root} of the edge, and each vertex is the root of a unique edge, called its {\it ball}. A De Bruijn covering code is
a subset of the roots such that every vertex is in some edge containing 
a chosen root. 
Under some mild conditions, we obtain an upper bound for the shortest length of a
De Bruijn covering code of a rooted hypergraph, a bound which is within a factor of $\log n$ of the lower bound.
\end{abstract}

\section{Introduction}

Suppose $G$ is a graph whose vertex set consists of some subset of all $n$-tuples $\cX^n$ over a finite alphabet $\cX$.  The natural distance metric $d(\cdot,\cdot)$ on this graph allows us, for each nonnegative integer $R$, to define a ``ball of radius $R$ centered at $x \in \cX$'' by
$$
B(x;R) = \{y : d(x,y) \leq R\}.
$$
One may ask for a subset of the vertices so that their respective balls cover the entire graph.  Such a set is commonly called a {\it covering code of radius $R$} for the graph $G$.

Because the vertex set of our graph consists of sequences of symbols, it is sometimes possible to find a particularly compact representation of a given covering code $\cC$.  Consider a string $S$ of elements of $\cX$ of length $|\cC|$.  We say that $S$ is a {\it De Bruijn covering code} for $G$ if the set of consecutive $n$-words (with ``wrap-around'') in $S$ is exactly the set $\cC$.  Then, instead of writing down all of $\cC$, we can specify the code with $n$ times the efficiency by using $S$.

Such an object was considered in \cite{ChC04} -- the graph was precisely the $q$-ary Hamming cube, $q$ a prime power, where our definition of a covering code coincides with the classical one.  In particular, the authors asked for the length of the shortest De Bruijn covering code of a given radius $R$, and showed that one exists with length given by
$$
\frac{q^n}{n^R} \prec |\cC| \prec \frac{q^n \log n}{n^R}.
$$
Later, Vu \cite{V05} extended these bounds to all $q$ and greatly simplified the proof.  The upper and lower bounds still do not meet, however, and it is an interesting question to close this gap.

It is also natural to ask analogous questions for other sets of sequences than all of $\cX^n$, possibly with an equivalence relation defined on these sequences.  One could ask for multisets of size $n$, permutations, or any Cayley graph $G$ defined on a subset of $\cX^n$.  For example, the string $134526$ is a radius $1$ covering code for the permutations on four symbols: indeed, every permutation is at a ``transposition distance'' of at most one from $1234$, $2341$, $2314$, $3241$, $2314$, or $4123$, the six order-types which occur as consecutive $4$-words in the string.  De Bruijn covering codes for these other types of graphs are precisely the subject of this paper.

Therefore we generalize this idea as follows.  A {\it hypergraph} is a pair 
$(\cV,\cE)$, where $\cV$ is the set of {\it vertices} and $\cE$ is a 
family of subsets of $\cV$, called {\it (hyper)edges}.  
A {\it $k$-uniform} hypergraph is one in which every 
hyperedge has cardinality $k$, and a {\it graph} is a $2$-uniform hypergraph. 
 Given a set $A$, 
a hypergraph $\cH$ on $A$ is said to be
a {\it rooted hypergraph} if every edge contains a special vertex, called the {\it root} of the edge, and each vertex is the root of exactly one edge.  Given a vertex $a \in A$, we denote the unique edge of which it is the root by $\hat{a}$, which we call the {\it ball about} $a$.  An {\it endomorphism} of a rooted hypergraph $\cH = (\cV,\cE)$ is a map $\phi : \cH \rightarrow \cH$ so that $\phi(e) \in \cE$ for all $e \in \cE$ and $\widehat{\phi(a)}=\phi(\hat{a})$ for all $a \in \cV$.  An {\it automorphism} of $\cH$ is a bijective endomorphism whose inverse is also an endomorphism.  $\cH$ is said to be {\it transitive} if its automorphism group acts transitively on its set of root-ball pairs.  We are primarily concerned with transitive rooted hypergraphs in the sequent.

A {\it covering code} for a rooted hypergraph $\cH$ on $A$ is a subset $S \subset A$ with the property that, for each $a \in A$, there exists a $b \in A$ with $a \in \hat{b}$.  In other words, a covering code is a set of vertices so that every vertex of the hypergraph belongs to some edge whose root lies in the set.

Suppose $\cX$ is a (finite or infinite) set, $\Pi$ is a family of disjoint subsets of $\cX^n$, and $\cH$ is a rooted hypergraph on $\Pi$.  Write $\Pi(\bx)$ for the member of $\Pi$ containing $\bx$.  Given a sequence $S = (s_0,\ldots,s_{M-1})$ with $s_i \in \cX$, write
$$
S^{(n)}_i = (s_{k},\ldots,s_{k+n-1}),
$$
and $S^{(n)} = \{S^{(n)}_i : 1 \leq i \leq M\}$, where all indices are taken modulo $M$.  We call $S$ an {\it order $n$ De Bruijn covering code} for $\Pi$ if $S^{(n)}$ is a covering code for $\cH$, and call $|S|=M$ its {\it length}.

We recover the previous definition of a De Bruijn covering code of radius $R$ by taking $\cX = \{0,1\}$, $\Pi$ the partition into singletons, and $\cH$ the set of Hamming $R$-balls, i.e., all $\hat{\bv} = \{\bw : d(\bv,\bw) \leq R\}$ for $\bv \in \{0,1\}^n$.\\

For a hypergraph $\cH$, we write $N(\pi)$ for the set of edges containing $\pi$, and $\deg(\pi)$ for $|N(\pi)|$.  Finally, we write $f(n) \prec g(n)$ to mean that there exists a $c > 0$ so that, for sufficiently large $n$, $f(n) \leq cg(n)$, and $f \sim g$ to mean that $g \prec f \prec g$.  Then we have the following theorem.

\begin{theorem} \label{mainthm} Let $(\cX,\Pi,\cH)$ be as above, with $\cX$ finite, and let $\cH$ be a transitive rooted hypergraph.  Suppose that $|\hat{\pi}| = K(n) \prec |\Pi|^2/n$ for all $\pi \in \Pi$.  Denote by $T_k$ the random variable $|\hat{\pi}_1 \cap \hat{\pi}_2|$, where $\pi_1$ and $\pi_2$ are chosen as follows: we pick a string $\Lambda$ uniformly at random from $\cX^{n+k}$, and set $\pi_1 = \Pi(\Lambda^{(n)}_1)$, $\pi_2 = \Pi(\Lambda^{(n)}_{k+1})$.  Then, if
$$
\sum_{1 \leq k \leq n-1} \bE(T_k) \prec K,
$$
there exists a De Bruijn covering code $S$ whose length satisfies
$$
\frac{|\Pi|}{K} \prec |S| \prec \frac{|\Pi| \log n}{K}.
$$
\end{theorem}

We delay the proof to Section \ref{mainproof}.  One useful corollary is the following simplification: take $\Pi$ to be the partition into singletons and the edges of $\cH$ to be the $R$-balls in the graphical distance metric.  The result then follows immediately from Theorem \ref{mainthm}, since $K \leq q^n \prec q^{2n}/n$ trivially.

\begin{cor} \label{graphthm} Let $G=(V,E)$ be a transitive graph, with $V=\cX^n$ for some set $\cX$ of cardinality $q < \infty$.  Suppose that $|B(\bv;R)| \sim K(n)$ for all $\bv \in V$.  Then, if
$$
\sum_{1 \leq k \leq n-1} \,\,\, \sum_{x_1,\ldots,x_{k+n} \in \cX} |B(x_1,\ldots,x_n;R) \cap B(x_{k+1},\ldots,x_{k+n};R)| \cdot q^{-n-k} \prec K,
$$
there exists a De Bruijn covering code $S$ whose length satisfies
$$
\frac{q^n}{K} \prec |S| \prec \frac{q^n \log n}{K}.
$$
\end{cor}

\section{De Bruijn covering codes for permutations and Hamming space}
Our first application generalizes the results of \cite{ChC04} and \cite{V05} to arbitrary (small) radii.

\begin{theorem} For $R = o(n)$, and any number of symbols $q$, there exists a De Bruijn covering code of radius $R$ for the $q$-ary Hamming space (with the ordinary Hamming metric) of dimension $n$ having cardinality $\prec q^n \log n / \binom{n}{R}$.
\end{theorem}

\begin{proof} We apply Corollary \ref{graphthm} to $\cX = \{1,\ldots,q\}$.  Clearly, for every $\bv \in \cX^n$, $|B(\bv;R)|=|B(\mathbf{0};R)|=\sum_{k=0}^R \binom{n}{k}$.  It is easy to see that
$$
\sum_{k=0}^R \binom{n}{k} \sim \binom{n}{R}
$$
when $R=o(n)$, so we may simply take $K = \binom{n}{R}$.  Now, fix $k$, $1 \leq k \leq n-1$, and consider the set
$$
S_k = \bigcup_{x_1,\ldots,x_{k+n} = 1}^q B(x_1,\ldots,x_n;R) \cap B(x_{k+1},\ldots,x_{k+n};R).
$$
We wish to count the elements of $S_k$.  They correspond to pairs of strings $(s,t)$, $s \in \cX^{n+k}$ and $t \in \cX^n$ so that $d(s_1^{(n)},t) \leq R$  and $d(s_{k+1}^{(n)},t) \leq R$.  In other words, the following family of equations holds:
\begin{align*}
s_1     =  t_1     =  s_{k+1} \\
\vdots \hspace{.5in} \\
s_n     =  t_n     =  s_{n+k}
\end{align*}
except for at most $R$ of the left-hand equalities and at most $R$ of the right-hand equalities.  For any particular choice of ``violated'' equalities, we may construct all possible solutions to the system by choosing $s_1,\ldots,s_k$ arbitrarily, and then arbitrarily choosing each of the values which follows an inequality.  All other values are determined by these choices.  That is, $|S_k| \prec \binom{n}{R}^2 q^{k+2R}$.

We therefore have that
\begin{align*}
\sum_{1 \leq k \leq n-1} \,\,\, |S_k| \cdot q^{-n-k} \prec \sum_{1 \leq k \leq n-1} \,\,\, \binom{n}{R}^2 q^{2R-n} = (n-1) \binom{n}{R}^2 q^{2R-n}.
\end{align*}
In order to apply Corollary \ref{graphthm}, we need $n \binom{n}{R} \prec q^{n-2R}$.  However, by Stirling's Formula,
$$
n \binom{n}{R} \prec \frac{n \cdot n^n}{R^R (n-R)^{(n-R)}} = n \cdot \left (  \epsilon^{-\epsilon} (1 - \epsilon)^{(\epsilon - 1)}  \right )^n,
$$
where $R = \epsilon n$.  Since $\lim_{\epsilon \rightarrow 0+} \epsilon^{-\epsilon} (1 - \epsilon)^{(\epsilon - 1)} = 1$, the result follows.
\end{proof}

We also have the following implication.  Given two sequences of reals $(a_1,\ldots,a_k) \in \mathbb{R}^k$ and $(b_1,\ldots,b_k) \in \mathbb{R}^k$, we say that they have the same {\it order type} if $a_i < a_j$ iff $b_i < b_j$ for every $1 \leq i,j \leq k$.  In \cite{ChDG92}, it is shown that there is a sequence of $n!$ reals so that the consecutive $n$-tuples represent every order type.  Here we have a similar result.

\begin{theorem} \label{firstapp} For any fixed positive integer $R$, there exists a sequence of reals $S = \{\alpha_0,\ldots,\alpha_{M-1}\}$ so that every permutation $\sigma$ on $n$ symbols differs from the order-type of some consecutive $n$-word in $S$ by at most $R$ transpositions, and
$$
\frac{n!}{n^{2R}} \prec M \prec \frac{n! \log n}{n^{2R}}.
$$
\end{theorem}
\begin{proof} Take $\cX = \mathbb{R}$, and $\Pi$ the equivalence relation on $n$-tuples of reals in $[0,1]$ {\it with no repeated elements} which represents their order-type  (i.e., $(x_1,\ldots,x_n) \sim_\Pi (y_1,\ldots,y_n)$ iff $x_i \leq x_j \Leftrightarrow y_i \leq y_j$ for all $i,j$.)  For $\pi \in P$, $\hat{\pi}$ is the $R$-ball rooted at the permutation $\pi$ under the transposition-distance metric.  Then
$$
K = |\hat{\pi}| = \frac{n^{2R}}{2^R R!} (1+o(1)) \prec |\Pi|^2/n = n!^2/n.
$$
In order to apply Theorem \ref{mainthm}, we must show that, if we choose $x_1,\ldots,x_{n+k}$ uniformly at random from $[0,1]$, then, if we denote by $\pi_1$ and $\pi_2$ the equivalence classes of $(x_1,\ldots,x_n)$ and $(x_{k+1},\ldots,x_{k+n})$, respectively, we have
$$
\sum_{1 \leq k \leq n-1} \bE(|\hat{\pi}_1 \cap \hat{\pi}_2|) \prec K.
$$
Fix $k \in [n-1]$.  If $|\hat{\pi}_1 \cap \hat{\pi}_2|$ is nonempty, then there is some set $S \subset [n]$ with $|S| \leq 4R$, so that whenever $\pi_1(x_1) < \pi_1(x_2)$ and $x_1,x_2 \in [n] \setminus S$, then $\pi_2(x_1) < \pi_2(x_2)$.  In particular, the event that $|\hat{\pi}_1 \cap \hat{\pi}_2| \neq \emptyset$ has probability at most $\binom{n}{2R}^2$ times the probability that two independent, uniformly chosen permutations of $[k] \setminus S$ are identical, since $\pi_1$ and $\pi_2$ restricted to this set are independent.  Therefore, for $n > 8R$,
$$
\sum_{8R \leq k \leq n-1} \bE(|\hat{\pi}_1 \cap \hat{\pi}_2|) \leq \binom{n}{2R}^3 \sum_{8R \leq k \leq n-1} \frac{1}{(k-4R)!} \leq \binom{n}{2R}^3 \frac{n}{(n/2)!} \prec K.
$$
Now, suppose $k < 8R$ and $n \geq 128R2 + 16R$.  Let $X_j = \{k(j-1)+1,\ldots,kj\}$ for $j = 1,\ldots,\floor{n/k}$.  Clearly, $\floor{n/k} \geq n/8R-1$.  On the other hand, only at most $4R$ of the $X_j$ contain a point of $S$.  Therefore, there is a run of consecutive $j$'s of length at least
$$
\frac{n/8R-1-k}{4R} \geq \frac{n - 8R - 64R2}{32R2} \geq \frac{n}{64R2}
$$
so that each $X_j$ contains no point of $S$.  This means that the least $x_i$ in each interval $\{x_{k(j-1)+1},\ldots,x_{kj}\}$ for $j$ in this range is a monotone sequence.  This event has probability $2/(n/64R2)!$, whence
$$
\sum_{1 \leq k < 8R} \bE(|\hat{\pi}_1 \cap \hat{\pi}_2|) \leq \frac{2}{(n/64R2)!} \binom{n}{2R} \prec K.
$$
\end{proof}

\section{Proof of the main theorem} \label{mainproof}

We need a result of Janson to proceed.  The following appears in \cite{Jan98}.  First, some notation.  Let $I$ be an index set for a set of events $\{B_i\}_{i \in I}$.  Define a graph $\sim$ on $I$ with the following property: Let $J_1$ and $J_2$ be two disjoint subsets of $I$ such that there is no $i_1 \in J_1$ and $i_2 \in J_2$ with $i_1 \sim i_2$.  Now, let $A_1$ be any Boolean function of the events $\{B_i : i \in J_1\}$ and let $A_2$ be any Boolean function of the events $\{B_i : i \in J_2\}$. Then $A_1$ and $A_2$ are independent.

Let $\mu = \sum_{i=1}^m \bP(B_i)$, $\Delta = \sum_{i \sim j} \bP(B_i \wedge B_j)$, and $\delta = \max_i \sum_{j \sim i} \bP(B_j)$.  Then the following holds.

\begin{lemma} With the above notation,
$$
\bP(\wedge_{i=1}^m \overline{B_i}) \leq \exp(-\min\left (\frac{\mu2}{8\Delta},\frac{\mu}{2},\frac{\mu}{6\delta}\right)).
$$
\end{lemma}

\begin{proof}[Proof of Theorem \ref{mainthm}.]  The number of vertices in $\cH$ is $|\Pi|$, so the lower bound follows immediately from the fact that $N$ edges cannot cover more than $O(KN)$ vertices.

As for the upper bound, our two-step strategy is as follows: first, we take a random set of vertices from the hypergraph; then, we ``patch up'' the string $S$ by appending all the equivalence classes we miss to the end.  So, take a string $S$ of length $M$, chosen from the uniform distribution on $\cX^M$.  Write $e_i$ for the element of $\Pi$ containing $S^{(n)}_i$, and let $\omega$ denote a random choice of $n$-word drawn uniformly at random.  Denote by $B_i^\pi$ the event that $\hat{e}_i$ contains $\pi \in \Pi$.  Then the probability that a given $\pi$ is covered by {\it no} edge whose root appears among the $e_i$ is given by $\bP(\wedge_{i=1}^m \overline{B_i^\pi})$.  Note that we may take $B_i^\pi \sim B_j^\pi$ iff $|i-j| \leq n \pmod{M}$.

We now estimate $\mu^\pi$, $\Delta^\pi$, and $\delta^\pi$.  Since $\deg(\pi) = K$ for each $\pi \in \Pi$ by the transitivity of $\cH$,
$$
\mu^\pi = \sum_{i=1}^M \bP(B^\pi_i) = \sum_{i=1}^M \bP(\pi \in \omega) = \sum_{i=1}^M \frac{\deg(\pi)}{|\Pi|} \sim \frac{MK}{|\Pi|},
$$
and
\begin{align*}
\Delta^\pi & = \sum_{1 \leq |i - j| \leq n \pmod{M}} \bP(B^\pi_i \wedge B^\pi_j) \\
&= \sum_{1 \leq |i - j| \leq n \pmod{M}} \frac{\bE(|\widehat{e}_i \cap \widehat{e}_j|)}{|\Pi|} \\
&= \sum_{i} \sum_{1 \leq k \leq n} \frac{\bE(T_k)}{|\Pi|} = \frac{M \bE(T_k)}{|\Pi|} \prec \frac{MK}{|\Pi|},
\end{align*}
so that
$$
\frac{\mu2}{8 \Delta} \prec \frac{MK}{|\Pi|}.
$$
Also,
$$
\delta = \max_i \sum_{j \sim i} \bP(B^\pi_j) \prec \frac{n \bE(|e_j|)}{|\Pi|} = \frac{n \sum_{\pi \in \Pi} |\hat{\pi}|}{|\Pi|^2} \sim \frac{nK}{|\Pi|^2} \prec 1.
$$
Therefore,
$$
\min\left (\frac{\mu2}{8\Delta},\frac{\mu}{2},\frac{\mu}{6\delta}\right) \gg \frac{MK}{|\Pi|},
$$
and setting $M = C |\Pi| \log n/K$ for an appropriate constant $C$ yields
$$
\bP(\wedge_{i=1}^m \overline{B_i^\pi}) \leq e^{-\log n} = n^{-1}.
$$
If $S$ is be a random string of length $M$, this means there are at most $|\Pi|/n$ classes in $\Pi$ which do not belong to any of the balls about points of $S^{(n)}$.  Then $S^\prime$ satisfies the conclusion of the theorem, where $S^\prime$ contains two consecutive copies of $S$ followed by a concatenated list of one string from each ``left out'' element of $\Pi$.
\end{proof}

\section{Problems and remarks}

The problem which looms largest is, of course, the elimination of the $\log n$ in the numerator of our upper bounds.  There are also a number of related questions which the techniques of this paper do not appear to resolve:

\begin{enumerate}
\item In \cite{ChDG92}, the authors show that a De Bruijn cycle exists for permutations of $n$ elements, using only $6n$ numbers.  That is, there is a sequence of length exactly $n!$, consisting of at most $6n$ reals, so that the set of order-types represented by consecutive $n$-words contains each permutation of $n$ exactly once.  We may ask something similar for De Bruijn covering codes: how long is the shortest De Bruijn covering code of radius $R$ for permutations of $n$ elements which uses only, say, $Cn$ symbols?

Furthermore, in \cite{ChDG92}, the authors conjectured there is a De Bruijn cycle for permutations of $n$ elements using exactly $n+1$ numbers.  This conjecture is still open.

\item An {\it error correcting} code, which can be thought of as dual to covering codes, is a subset of $[q]^n$ so that no two words in the code are less than $R$ symbol-changes apart.  What is the longest $q$-ary sequence whose consecutive $n$-words form an error correcting code of radius $R$?  Clearly, we may ask analogous questions for permutations and other rooted hypergraphs.
\item Given a subset $S$ of the permutations on $n$ symbols and a radius $R$, one may ask, what is the length of the shortest De Bruijn covering code for $S$?  If $|S|=o(n)$, then choosing the sequence uniformly at random may be very suboptimal.  For example, one may take $S$ to be the set of permutations with at most $k$ descents.
\end{enumerate}

\end{document}